\documentclass[12pt]{article}
\usepackage{epsfig,color,oldgerm,amsmath,pifont,amsfonts,amssymb,latexsym,umlaut}
\usepackage{float,xspace,calc,theorem,ifthen,enumerate,a4wide}

\usepackage{hyperref}
{\theoremstyle{break} }
\newtheorem {lemma}{Lemma} 
\newtheorem{prop}{Proposition} 
\newtheorem {koro}{Corollary} 
\newtheorem{theo}{Theorem} 
{\theorembodyfont{\rmfamily}\theoremstyle{plain}
 
}
{\theorembodyfont{\rmfamily}\theoremstyle{break}

\newcommand{\qed}{\phantom{xxxxxx}\hfill q.e.d.}
\newenvironment{proof}{{\noindent\em Proof:} }{\qed\\}
\newenvironment{proofof}[1]{{\noindent\em Proof of #1:} }{\qed\\}

\newcommand{\ess}{{\rm ess }} 
\newcommand{\bfone}{{\bf 1}} 
\newcommand{\const}{{\rm const\,}} 
\newcommand{\var}{\operatorname{var}} 
\newcommand{\sign}{\operatorname{sign}} 
\newcommand{\ph}{\varphi}
\renewcommand{\th}{\vartheta}
\newcommand{\esup}{\operatornamewithlimits{ess\,sup}}
\newcommand{\einf}{\operatornamewithlimits{ess\,inf}}
\newcommand{\bvnorm}[1]{\|#1\|_{BV}}
\newcommand{\BV}{\mathit{BV}}

\newcommand{\rz}{{\mathbb R}}
\newcommand{\nz}{{\mathbb N}}
\newcommand{\cz}{{\mathbb C}}
\newcommand{\zz}{{\mathbb Z}}

\newcommand{\tz}{{\mathbb T}}
\renewcommand{\SS}{\cal S}

\newcommand{\eg}{{e.g.}\xspace} 
\newcommand{\ie}{{i.e.}\xspace}

\newcommand{\rr}{\rightarrow}

\begin{document}
\title{Eigenfunctions for smooth expanding circle maps}
\author{Gerhard Keller\thanks{Universit\"at Erlangen-N\"urnberg,
 Germany, {\tt keller@mi.uni-erlangen.de}}\, 
and Hans-Henrik Rugh\thanks{Universit\'e de Cergy-Pontoise,
{}France, {\tt hhrugh@math.u-cergy.fr}}
\hspace*{3mm}\thanks{G.K. acknowledges the support and hospitality of the Mathematics Department of the 
Universit\'e de Paris-Sud, France, where part of this work was done.}}
\maketitle

\begin{abstract} We construct a real-analytic circle map for which
  the corresponding Perron-Frobenius operator has a real-analytic
  eigenfunction with an eigenvalue outside the essential spectral
  radius when acting upon $C^1$-functions.
\end{abstract}

\section{Introduction}
\label{sec:intro}
Let $T$ be a uniformly expanding $C^2$ map of the
circle $S^1=\rz/\zz$. Consider
an associated Perron Frobenius operator acting upon $C^1$ functions,
\begin{equation}
    {P}f(x) = \sum_{y: Ty=x} \frac{1}{|T'(y)|} f(y), \ \ x\in S^1.
\end{equation}
By the Ruelle-Perron-Frobenius Theorem, one is a simple eigenvalue and
the operator has a `spectral gap'. It is also known \cite{CI91} (see also
\cite[Th.~2.5]{baladi}) that the essential spectrum of ${P}$ is a disk centered at zero and of radius $r_{\ess}=\th$ with
$\th^{-1}=\liminf_k \inf_x \sqrt[k]{|(T^k)'(x)|}$.  Outside this
disk there can only be isolated eigenvalues of finite
multiplicity. A natural question is then if there {\em really are} any
eigenvalues apart from one outside the essential spectrum?

This `non-essential' spectrum is void in typical standard examples.
{}For example, $Tx=2x \mod 1$, with ${P}$ acting upon real-analytic
functions has 1 as a simple eigenvalue and the rest of the spectrum
concentrated at zero. On $C^1$ functions the essential spectral radius
is one-half and a small perturbation can thus not give rise to an
example with a non-trivial non-essential spectrum. Even if we consider
$Tx = 2 x \mod 1$ with ${P}$ acting on $C^1$ function on the interval
$[0,1]$ the largest eigenvalue below 1 is one-half, precisely on the
boundary of the essential spectrum. More generally, if
$T$ is a smooth
expanding circle map fixing say zero we may for $k\geq 1$ consider
the Perron Frobenius operator as an operator $P_{S^1}$ acting upon
 $C^k(S^1)$, 
or as $P_{[0,1]}$ acting upon $C^k([0,1])$.
The spectrum of $P_{[0,1]}$ is the same as the spectrum of $P_{S^1}$
union a point spectrum $\{|T'(0)|^{-m}, m\geq 1\}$. In particular,
in the case $k=1$ the spectrum outside the
radius $\th$ doesn't change.
One may show this by adapting the proof of 
Corollary 2.11 in \cite{Rue}.

  Recently, Liverani \cite{Liv} suggested to us that when
acting upon $C^1$ functions there may not be any non-trivial
`non-essential' spectrum.  The present paper was conceived in an
attempt to resolve this question and we have arrived at the following

  \begin{theo}\label{theo:main} 
      We may construct a  real-analytic circle map
 $T: S^1\rr S^1$ of derivative greater than $3/2$ for which ${P}$
  has a real-analytic eigenvector with eigenvalue 
  greater than $0.75$ in absolute value.
   In particular, this eigenvalue is greater than the
  essential spectral radius for $C^1(S^1)$ functions 
 which does not exceed $2/3$.
  \end{theo}

Our proof consists 
of a fairly explicit example homotopic to
the $2x \mod 1$ map. We will construct it in two steps,
(1) we will find
a piecewise linear map acting upon functions of bounded variation
($\BV$) and (2) we apply a real-analytic smoothening of this map. 
The advantage when acting upon $\BV$ is that this space allows for 
discontinuities in the derivatives of the dynamical system and in
the eigenfunctions of the associated operator. This makes it very easy 
to construct examples with very explicit spectral properties.
The hard part is to recover a corresponding
 result in the smooth category. We do this  through what may be called
`quasi-compact' perturbation theory for $\BV$ operators.

\textbf{Acknowledgment:} We thank Viviane Baladi and Carlangelo
Liverani for bringing this problem to our attention and for
stimulating discussions about it.

\section{Rychlik's result}
\label{sec:rychlik}

Although there are numerous treatments of piecewise expanding interval
maps in various contexts, Rychlik's early paper \cite{rychlik} is
still among the most elegant and general ones and seems most suitable
as a base for our arguments. We recall his setting.  Let $X$ be the
closed unit interval, $m$ Lebesgue measure on $X$, $U\subseteq X$ open
and dense with $m(U)=1$ and denote $S=X\setminus U$. We assume that
$S$ contains the endpoints of $X$ so that $U$ is an at most countable
union of open intervals. The closures of these intervals are then a
countable family, call it $\beta$, of closed intervals with disjoint
interiors such that $\bigcup_{B\in\beta}\beta\supset U$ and such that
for any $B\in\beta$ the set $B\cap S$ contains exactly the endpoints
of $B$.

The space of functions of bounded variation $f:X\rr\rz$ is denoted by
$\BV$, the variation of $f$ over all of $X$ by $\var(f)$ and the
variation over a subinterval $I\subseteq X$ by $\var_I(f)$.

Let $T:U\rr X$ be a continuous map such that for any $B\in\beta$ the
restriction $T|_{B\cap U}$ admits an extension to a homeomorphism of
$B$ with some interval in $X$.

Suppose that $T|_U$ is differentiable in the following sense: There is
a function $g:X\rr\rz_+$ with $\|g\|_\infty<\infty$, $\var(g)<\infty$
and $g|_S=0$ such that the transfer operator $P$ defined by
\begin{displaymath}
  Pf(x)=\sum_{y\in T^{-1}x} g(y)f(y)
\end{displaymath}
preserves $m$ which means $m(Pf)=m(f)$ for all bounded measurable
$f:X\rr\rz$.

{}For iterates $T^N$ of $T$ we adopt the following notation:
$S_N=\bigcup_{k=0}^{N-1}T^{-k}S$, $U_N=X\setminus S_N$ and
$\beta_N=\bigvee_{k=0}^{N-1}T^{-k}\beta$. Next define $g_N:X\rr\rz_+$
by $g_N|_{S_N}=0$ and $g_N|_{U_N}=g\circ T^{N-1}\cdot\dots\cdot g\circ
T\cdot g$. In order to see that $T^N,U_N,S_N,g_N,\beta_N$ satisfy our conditions for $N\geq1$ the only nontrivial thing to verify is $\var(g_N)<\infty$:
\begin{lemma}\protect{\cite[Lemma 2]{rychlik}}
  \begin{displaymath}
    \var(g_N)\leq 2^{N-1}(\var g)^N. 
  \end{displaymath}
  (Indeed, since we assumed that the endpoints of $X$ belong to $S$,
  inspection of the proof reveals that $\var(g_N)\leq(\var g)^N$.)
\end{lemma}

\begin{lemma}\protect{\cite[Lemma 6]{rychlik}}
  \label{lemma:Ry6}
  {}For every $\epsilon>0$ and every $N\geq1$ there exists a finite partition 
  $\alpha_N$ of $X$ into intervals such that
  \begin{displaymath}
    \max_{A\in\alpha_N}\var_A(g_N)<\|g_N\|_\infty+\epsilon\ .
  \end{displaymath}
\end{lemma}
\begin{koro}
\label{koro:Ry6}
  In the situation of the previous lemma there exists
  $\rho=\rho(\alpha_N)>0$ such that actually
  $\var_{U_\rho(A)}(g_N)<\|g_N\|_\infty+\epsilon$ where $U_\rho(A)$
  denotes the $\rho$-neighbourhood of $A$.
\end{koro}
\begin{lemma}\protect{\cite[Corollary 3 and its proof]{rychlik}}
  \label{lemma:RyCo3}
  {}For every $N\geq1$ and $\lambda_N>2\|g_N\|_\infty$ there exists
  $D_N\geq0$ such that for every $f:X\rr\rz$
  \begin{displaymath}
    \var(P^Nf)
    \leq
    \lambda_N\var(f)+D_N\|f\|_1
  \end{displaymath}
  where $D_N$ is determined as follows: By Lemma~\ref{lemma:Ry6} there
  is a finite partition $\alpha_N$ such that
  $\max_{A\in\alpha_N}\var_A(g_N)\leq\lambda_N-\|g_N\|_\infty$. Then
  $D_N=\max_{A\in\alpha_N}\var_A(g_N)/m(A)$.
\end{lemma}

{}From now on we assume that $\|g_N\|_\infty<1$ for some $N\geq1$.  Then
\begin{displaymath}
  \th:=\lim_{n\rr\infty}\|g_N\|_\infty^{1/N}<1.    
\end{displaymath}

\begin{lemma}\protect{\cite[Proposition 1 and its proof]{rychlik}}
\label{lemma:RyProp1}
  Given $\kappa\in(\th,1)$ we can find $F\geq0$ such that for every
  $f:X\rr\rz$ and every $n\in\nz$
  \begin{displaymath}
    \sum_{B\in\beta^n}\var P^n(f\cdot\chi_B)
    \leq
    F\cdot(\kappa^n\var(f)+\|f\|_1)\ .
  \end{displaymath}
  $F$ is determined as follows: Fix $M$ such that
  $\th<\|2g_M\|_\infty^{1/M}<\kappa<1$. Then
  $F:=\max\{D/(1-\kappa^M),\lambda/\kappa^{M-1}\}$ where
  $\lambda:=\max\{\lambda_1,\dots,\lambda_M\}$ and
  $D:=\max\{D_1,\dots,D_M\}$ with $\lambda_i,D_i$ from
  Lemma~\ref{lemma:RyCo3}.
\end{lemma}
\begin{koro}\label{koro:LY}
  With $F$ and $\kappa$ as before we have for every $f\in\BV$
  \begin{displaymath}
    \var(P^nf)
    \leq
    F\cdot(\kappa^n\var(f)+\|f\|_1)\ .
  \end{displaymath}
\end{koro}

It is this explicit knowledge of the dependence of $F$ and $\kappa$ on $T$
and $g$ which will allow us later to apply spectral perturbation theory to smooth
approximations of piecewise linear maps.

\section{Smoothing}
\label{sec:smoothing}
We want to apply the Lasota-Yorke type estimate of
Corollary~\ref{koro:LY} to piecewise linear circle maps and their
smooth approximations. In terms of interval maps this will be full
branched maps. To this end fix $p\in\zz$, $p\geq2$, and let
%HHR 
$\psi: \rz/p\zz    \rr   \rz/\zz$
be an increasing homeomorphism. It lifts to a continuous strictly 
increasing  map
$\tau:\rz\rr\rz$ for which
$\tau(x+p)=\tau(x)+1$ for all $x$.  Let $S=\tau(\zz\cap{[0,p]})$,
$U=X\setminus S$, and denote by $T:U\rr X$ the $p$-branched map
defined by
\begin{displaymath}
  T(x)=\tau^{-1}(x)\mod 1
\end{displaymath}
(Of course $T$ can be interpreted as a continuous $p$-fold covering
map of $\tz^1$.) 

In what follows we assume that $\tau$ is twice differentiable in the
following sense: There is a function $\dot\tau:\rz\rr\rz$ of locally
bounded variation and with period $p$ such that
$\tau(x)=\int_0^x\dot\tau(u)\,du$ for all $x$. We assume that
$0<\inf\dot\tau\leq\sup\dot\tau\leq1$, and we set
$g|_U:=\dot\tau\circ\tau^{-1}|_U$ and $g|_S=0$. Then $T$, $g$ and
$\beta:=\{[\tau(k-1),\tau(k)]:k=1,\dots,p\}$ fit the setting of
section~\ref{sec:rychlik}.

Let $\ph:\rz\rr[0,\infty)$ be a smooth convolution kernel, \ie
$\lim_{|x|\rr\infty}\ph(x)=0$, $\int_{-\infty}^\infty\ph(x)\,dx=1$.
To be definit we choose $\ph$ to be the density of the standard
normal distribution, in which case $\ph$ is indeed real-analytic.  For
$\delta>0$ define the rescaled kernels $\ph_\delta$ by
$\ph_\delta(x)=\delta^{-1}\ph(\delta^{-1}x)$. We use these kernels to
define smooth approximations
$\tau_\delta:=\tau*\ph_\delta-\tau*\ph_\delta(0)$ to $\tau$ and
$\dot\tau_\delta:=\dot\tau*\ph_\delta$ to $\dot\tau$. Obviously
$\tau_\delta(0)=0$, $\tau_\delta(x+p)=\tau_\delta(x)+1$, and
$\dot\tau_\delta(x+1)=\dot\tau_\delta$. A simple calculation yields
\footnote{ $\int_0^x\dot\tau_\delta(u)\,du =
  \int_0^x\int_{-\infty}^\infty\dot\tau(y)\ph_\delta(u-y)\,dy\,du =
  \int_0^x\int_{-\infty}^\infty\tau(y)\ph_\delta'(u-y)\,dy\,du=\\ =
  \int_{-\infty}^\infty\int_0^x\tau(y)\ph_\delta'(u-y)\,du\,dy=
  \int_{-\infty}^\infty\tau(y)\,(\ph_\delta(x-y)-\ph_\delta(-y))\,dy=
  \tau*\ph_\delta(x)-\tau*\ph_\delta(0)=\tau_\delta(x)$.  }
\begin{equation}
\label{eq:tau_delta}
  \tau_\delta(x)=\int_0^x\dot\tau_\delta(u)\,du
\end{equation}

Denote by $T_\delta$ the $p$-branched transformation of $X$
determined by $\tau_\delta$ and define $g_\delta$ in terms of
$\tau_\delta$ as $g$ was defined in terms of $\tau$. Denote by
$P_\delta$ the transfer operator for $T_\delta$ and $g_\delta$.

Observe that 
\begin{displaymath}
  \label{eq:bounds}
  0<\einf\dot\tau\leq\tau_\delta'\leq\esup\dot\tau\leq1
\end{displaymath}
so that also
\begin{displaymath}
  0<\einf g\leq g_\delta\leq\esup g\leq1\ .
\end{displaymath}
In particular 
\begin{equation}
\label{eq:g_M}
  \|2g_{\delta,M}\|_\infty^{1/M}
  \leq
  2^{1/M}\cdot\|g_M\|_\infty^{1/M}
\end{equation}
so that the integer $M$ in Lemma~\ref{lemma:RyProp1} can be chosen
uniformly for $\delta\in(0,1)$ provided $\kappa\in(\th,1)$.
Then all $\lambda_N$ in Lemma~\ref{lemma:RyCo3} can be chosen to be
equal to $3$ so that also $\lambda=3$ in Lemma~\ref{lemma:RyProp1}. In
view of Corollary~\ref{koro:Ry6} the same partitions $\alpha_N$ can be
used in Lemma~\ref{lemma:RyCo3} for all $T_\delta$, $g_\delta$
provided $\delta>0$ is sufficiently small ($N=1,\dots,M$). Hence also
the numbers $D_N$ can be chosen uniformly for sufficiently small
$\delta$ ($N=1,\dots,M$). This proves:
\begin{prop}
  \label{prop:LY}
  The constants $F$ and $\kappa$ in the Lasota-Yorke type inequality
  \begin{displaymath}
    \var(P_\delta^nf)
    \leq
    F\cdot(\kappa^n\var(f)+\|f\|_1)\ .
  \end{displaymath}
  for $P_\delta$ (compare to Corollary~\ref{koro:LY}) can be chosen
  uniformly in $\delta$ provided $\delta>0$ is small enough.
\end{prop}

This Proposition opens the way to apply the spectral perturbation
theorem of \cite{KL99} to $P$ and $P_\delta$. We just need to show
that, for sufficiently small $\delta>0$, the operator $P_\delta$ is
close to $P$ in a suitable sense. The following proposition provides
the relevant estimate.
\begin{prop}\label{prop:approximation}
  {}For all $\delta\in(0,1)$ and
  all $f\in\BV$:
  \begin{displaymath}
    \int|P_\delta f(x)-Pf(x)|\,dx\leq 2C_0\cdot\delta\cdot\bvnorm{f}
  \end{displaymath}
  where $\bvnorm{f}:=\var(f)+\|f\|_1$.
\end{prop}
{}For its proof we need some simple estimates.
\begin{lemma}
  There is a constant $C_0>0$ (depending on $\tau$ and on the kernel
  $\ph$) such that for all $\delta\in(0,1)$:
\begin{enumerate}[a)]
\item
  $\sup_{x\in{}X}|\tau_\delta^{-1}(x)-\tau^{-1}(x)|\leq{}C_0\cdot\delta$
\item $\int_0^p|\tau_\delta'(y)-\dot\tau(y)|\,dy\leq{}C_0\cdot\delta$
\end{enumerate}
\end{lemma}
\begin{proof}
  \begin{enumerate}[a)]
  \item As $(\einf\dot\tau)^{-1}$ and $\esup\dot\tau$ are Lipschitz
    constants for $\tau^{-1}$ and $\tau$ respectively, we have
    \begin{align*}
      \sup_{x\in X}|\tau_\delta^{-1}(x)-\tau^{-1}(x)|
      &=
      \sup_{y\in[0,p]}|y-\tau^{-1}(\tau_\delta y)|
      \leq
      \frac1{\einf\dot\tau}\cdot\sup_{y\in[0,p]}|\tau(y)-\tau_\delta(y)|\\
      &\leq
      \const\cdot\delta\cdot\frac{\esup\dot\tau}{\einf\dot\tau}
      =:C_0\cdot\delta
    \end{align*}
  \item This follows from the assumption that $\dot\tau$ is of bounded
    variation.
  \end{enumerate}
\end{proof}
\begin{proofof}{Proposition~\ref{prop:approximation}}
  The proof is similar to the one of Lemma 13 in \cite{keller82}.
  Let $\psi:=\sign(P_\delta f-Pf)$. Then 
  \begin{equation}\label{eq:diff1}
    \begin{split}
      \int_X|P_\delta f(x)-P f(x)|\,dx 
      &= 
      \int_X(P_\delta f(x)-Pf(x))\,\psi(x)\,dx \\
      &=
      \int_Xf\cdot(\tilde\psi(\tau_\delta^{-1}x)-\tilde\psi(\tau^{-1}x))\,dx
    \end{split}
  \end{equation}
  where $\tilde\psi$ is the periodic extension of $\psi$ from $X$ to
  all of $\rz$. 
  
  Observe now that for each $\ph:X\rr\rz$ which is differentiable in
  the sense that $\ph(x)=\int_0^x\dot\ph(u)\,du$ for some bounded
  measurable $\dot\ph:X\rr\rz$ we have $\int_X
  f(x)\dot\ph(x)\,dx\leq\var(f)\,\sup_X|\ph|$.
  Therefore we can continue (\ref{eq:diff1}) by
  \begin{align*}
    \int_X|P_\delta f(x)-Pf(x)|\,dx 
    &\leq
    \bvnorm{f}\cdot\sup_{x\in X}\left|\int_0^x\tilde\psi(\tau_\delta^{-1}u)-\tilde\psi(\tau^{-1}u)\,du\right|\\
    &=
    \bvnorm{f}\cdot\sup_{x\in X}\left|\int_0^{\tau_\delta^{-1}x}
      \tilde\psi(y)\tau_\delta'(y)\,dy
      -\int_0^{\tau^{-1}x}\tilde\psi(y)\dot\tau(y)\,dy\right|\\
    &\leq
    \bvnorm{f}\cdot\left(\sup_{x\in X}|\tau_\delta^{-1}x-\tau^{-1}x|
      +\int_0^p|\tau_\delta'(y)-\dot\tau(y)|\,dy\right)\\
    &\leq
    2C_0\cdot\delta\cdot\bvnorm{f}
  \end{align*}
\end{proofof}

\section{A piecewise linear example}
\label{sec:pwl}

\begin{figure}
  \begin{center}
  \epsfig{figure=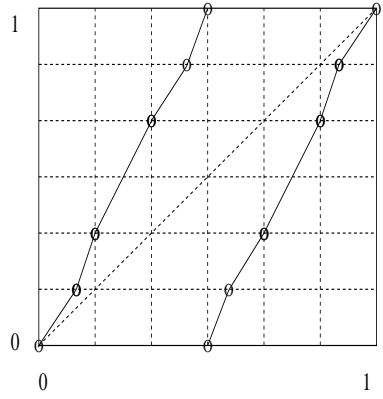,width=6cm,height=6cm}
   \caption{The piecewise linear map $T$}
  \label{fig:graph}
   \end{center}
\end{figure}

Let now $p=2$ and denote by $\psi: \rz/2\zz \rr \rz/\zz$ the
homeomorphism fixing zero and having successive slopes
$\frac23,\frac13, \frac12,\frac12, \frac23,\frac13, \frac13,\frac23,
\frac12,\frac12, \frac13,\frac23 $ on the intervals
$(\frac{k-1}6,\frac k6) ({\rm mod}\ 2)$, $k=1,\dots,12$.  The lift
$\tau: \rz \rr \rz$ of $\psi$ satisfies $\tau(x+2)=\tau(x)+1$.  Also
let $\pi : \rz/2\zz \rr \rz/\zz$ be the canonical projection.  The
map $T = \pi \circ \psi^{-1} : \rz/\zz \rr \rz/\zz$ is a piecewise
linear uniformly expanding map of the circle (cf. Figure \ref{fig:graph}).

Denote $I_k=(\frac{k-1}6, \frac k6)\ {\rm mod}\ \zz$, $k=1,\ldots,6$ 
and $C=\frac16 \zz/\zz$.
We write $S=\psi(\frac16 \zz/2\zz)$ for the 12 points on the circle where
$(T'(x))^{-1}=\psi'\circ \psi^{-1}(x)$ is not defined. We set
\begin{equation}
   g(x) = \left\{ 
        \begin{array}{cl}
           \psi' \circ \psi^{-1} (x)  \  & x\in \rz\setminus S, \\
            0                          & x \in S .
         \end{array}
        \right.
\end{equation}
As before we have a corresponding Perron Frobenius operator
\begin{equation}
P f(x) = \sum_{y:Ty=x} g(y) f(y) .
\end{equation}
 $P$ maps $\BV$ into the subspace
\begin{equation}
     E_C = \{ \phi \in \BV : \phi_{|C} \equiv 0\}.
\end{equation}

$P$ also preserves the subspace of step functions
\begin{equation}
     L_C = \{ \phi = \sum_{i=1}^6 c_i \bfone_{I_k} 
                : c = (c_i)_{i=1..6} \in \rz^6\}\ :
\end{equation}
acting with the operator $P$ upon $L_C$ we obtain $P\sum_i c_i
\bfone_{I_i} = \sum_{ij} c_i M_{ij} \bfone_j$ where $M$ is the doubly
stochastic matrix
\begin{displaymath}
  M = \begin{pmatrix}
    2/3&1/3&0&0&0&0\\
    0&0&1/2&1/2&0&0\\
    0&0&0&0&2/3&1/3\\
    1/3&2/3&0&0&0&0\\
    0&0&1/2&1/2&0&0\\
    0&0&0&0&1/3&2/3
  \end{pmatrix}
\end{displaymath}
Therefore, if $v$ is a left eigenvector of this matrix,
$\lambda{}v=vM$, then it induces an eigenfunction $\phi_v$ of
$P:\BV\rr\BV$ with the same eigenvalue $\lambda$. More precisely,
$\phi_v=\sum_i \mu_i\,\bfone_{I_i}$. In particular, $P1=1$.

The matrix $M$ has simple eigenvalues $\lambda_1=1$,
$\lambda_2=-\frac16-\frac{\sqrt{13}}{6}\approx-0.7676$,
$\lambda_3=\frac23$,
$\lambda_4=-\frac16+\frac{\sqrt{13}}6\approx0.4343$ and a double
eigenvalue at $0$.  As $0$ is a fixed point with slope $\frac32$, the
essential spectral radius of $P$ is $\th=\frac23$.  It follows that
the spectrum of $P:\BV\rr\BV$ is contained in the set
$\{1,\lambda_2\}\cup\{|z|\leq\th\}$, see \eg \cite[Lemma 3.1]{BK97}.
Indeed, a straightforward generalization of the the proof of that
lemma shows that $\lambda_2$ is a simple eigenvalue of $P$ because it
is a simple eigenvalue of $M$ of modulus larger than $\th$.  
Let $\Phi_2:\BV\rr\BV$ be the corresponding spectral projector.

The unique (up to normalization) left eigenvector of $M$ with
eigenvalue $\lambda_2$ is
$v_2=(1,\frac{3+\sqrt{13}}2,-\frac{5+\sqrt{13}}2,-\frac{5+\sqrt{13}}2,\frac{3+\sqrt{13}}2,1)$.
The associated eigenfunction $\phi_{v_2}$ is piecewise constant and
its set of essential discontinuities is
$D=\{\frac16,\frac26,\frac46,\frac56\}$.  Observe that these are all
essential discontinuities also if $\phi_{v_2}$ is considered as a
function from $S^1$ to $\rz$. (The discontinuity at $\frac 12$, where
$\phi_{v_2}(\frac12)=0$ by definition, is inessential, because the
left and right limits of $\phi_{v_2}$ at this point are both
$-\frac{5+\sqrt{13}}2$. Similarly the discontinuity at $0\in S^1$ is
inessential.)

\begin{lemma}\label{lemma:non-zero}
  There is a trigonometric polynomial $f:X\rr\rz$ such that
  $\Phi_2(f)\neq0$.
\end{lemma}
\begin{proof}
The difficulty arises from the 
fact that smooth functions are not
dense in $BV$. 
Thus a non-zero linear functional on $BV$ may vanish when
acting upon any smooth
function. We need to show that this is not the case for 
the eigenprojection $\Phi_2$.

Consider the space $E_C^\infty =
 \{ \phi \cdot \bfone_{S^1\setminus C} : \phi\in C^\infty(S^1) \}$.
We claim that 
$\phi_{v_2} \in E_C^\infty + P E_C^\infty$. To see this let 
$\psi\in C^\infty(S^1)$ be a test function with support in a 
`sufficiently small'
neighborhood of $0$. 
{}For example, we may suppose that $\psi(x)\equiv 1$
when $d(x,0)\leq1/100$ and that 
$\psi(x)=0$ for $d(x,0)\geq 1/40$.
{}For $x\in S^1$ we write
$\psi_x(y)=\psi(x+y)$ for the translated test-function.

\begin{figure}
  {\input{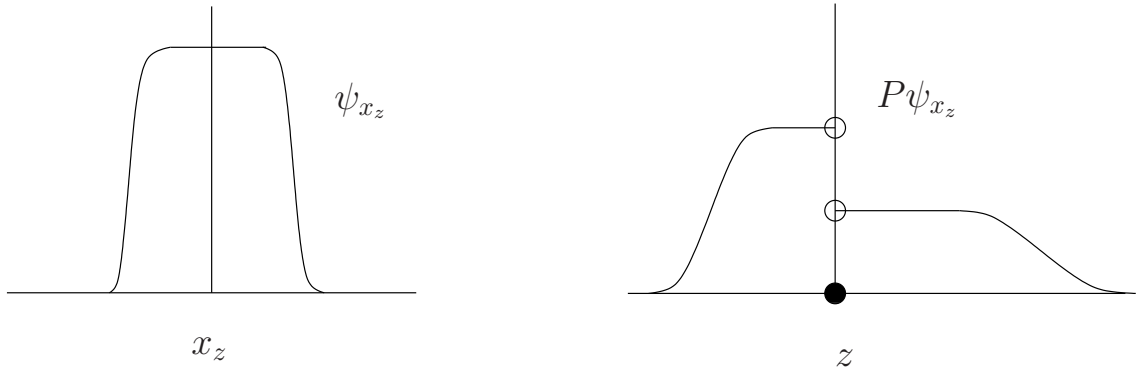}}
  \caption{A test-function and its image under $P$}
  \label{fig:2}
\end{figure}

Let $D$ be the above  set of essential discontinuities of 
 $\phi_{v_2}$. For each $z\in D$ we denote by $x_z$ the unique
T-preimage of $z$ in $(0,\frac12)$. 
Consider now the action of $P$ upon our test-function $\psi_{x_z}$.
The support of $\psi_{x_z}$ is so small (this is the meaning of
`sufficiently small' above) that
$P\psi_{x_z}$ has an essential discontinuity at $z$ but is otherwise
$C^\infty$, see Figure~\ref{fig:2}.
Since $\phi_{v_2}$ is locally constant on each side of $z$
we may find $\alpha_z\neq 0$ such that
$\phi_{v_2}-\alpha_z P\psi_{x_z}$ is locally constant in a 
punctured neighborhood of $z$ (simply let $\alpha_z$ be the ratio
of essential discontinuities at $z$ of the two functions).
Doing this for all $z\in D$ we see that
  \begin{equation}
       \tilde{\phi} = \phi_{v_2} - \sum_{z\in D} \alpha_z P\psi_{x_z}
  \end{equation}
  is locally constant at each punctured neighborhood of $z\in C$,
  vanishes at $z\in C$, and is $C^\infty$ elsewhere. In other words
  $\tilde{\phi}\in E_C^{\infty}$ and we have obtained the desired
  decomposition
  \begin{equation}
     \phi_{v_2}=\tilde{\phi} + P \phi, \ \ \tilde{\phi},\phi \in E_C^\infty.
  \end{equation}

Assume now that
$C^\infty(S^1)$ was in the kernel of $\Phi_2$. Functions with support in $C$
are in the kernel of $P$, and therefore also of $\Phi_2$. But then every
$\phi\in E_C^\infty$ is in the kernel of $\Phi_2$. The calculation,
 \begin{equation}
       \phi_{v_2} = \Phi_2 \phi_{v_2} =
       \Phi_2 \tilde{\phi} + \Phi_2 P \phi =
       0 + P \Phi_2 \phi = 0
  \end{equation}
is a contradiction, showing that $\Phi_2$ can not vanish upon
the space of $C^\infty$ functions.
  Since $C^\infty$ functions can be
  approximated in $\BV$-norm by trigonometric polynomials, there is
  also a trigonometric polynomial $f$ with $\Phi_2(f)\neq0$.
\end{proof}

\section{ Proof of Theorem~\ref{theo:main}}
\label{sec:proof}

We will construct the example announced in Theorem~\ref{theo:main} as
a smooth perturbation of the above piecewise linear map $T$.  Indeed,
for $\th=\frac23<\kappa:=0.7<\lambda_{2}$ we can apply the reasoning
of section~\ref{sec:smoothing} to $T$. The spectral perturbation
theorem in \cite{KL99} then shows that for sufficiently small
$\delta>0$ the transfer operator $P_\delta$ of the smooth map
$T_\delta$ (as an operator on $\BV$) has essential spectral radius at
most $\kappa$ and has exactly two simple eigenvalues $1$ and
$\lambda_\delta\equiv \lambda_{2,\delta}$ close to $\lambda_2$ outside
the essential spectrum.  The associated eigenprojector
$\Phi_{2,\delta}$ is a small perturbation (in $\BV$-norm) of the
eigenprojector $\Phi_2$. In particular it
does not vanish on the space of trigonometric polynomials if
$\delta>0$ is sufficiently small, see Lemma~\ref{lemma:non-zero}. We
want to show that associated to $\lambda_{\delta}$ there is an
eigenfunction of $P_\delta$ which is real-analytic. This will
prove our Theorem.\\

{}First note that there is at least one left eigenvector, $\ell_\delta\in
\BV^*$ associated to $\lambda_\delta$, i.e.\ $\ell_{\delta} P_\delta =
\lambda_{\delta} \ell_{\delta}$.  (One may take
$\ell_\delta(f)=\int\Phi_{2,\delta}(f)h_{2,\delta}\,dm$ where
$h_{2,\delta}$ is the eigenfunction of $P_\delta$ with eigenvalue
$\lambda_{2,\delta}$.) Since $\Phi_{2,\delta}$ and $h_{2,\delta}$ are
small perturbations of the corresponding objects $\Phi_2$ and
$h_2=\phi_{v_2}$, we may choose $\delta>0$ so small that $\ell_\delta$
does not vanish on the space of trigonometric polynomials on $S^1$.
We fix from now on a value of $\delta>0$ for which the above holds.

The map $\tau_\delta$, smoothened by convolution with the Gaussian
kernel, is real-analytic and has derivative smaller than $2/3$
(because, a.e.\ the derivative of $\psi$ varies between $1/3$ and
$2/3$).  We may therefore find $\rho>0$ and $\th<1$ (close to $2/3$)
for which $|\tau_\delta'(z)| \leq \th$ for all $z$ in the annulus
$A^p_\rho = \{ z \in \cz/p \zz : | {\rm Im}\,z | < \rho \}$.  Then
$\tau_\delta$ is a contraction from $A^2_\rho$ into $A^1_{\th \rho}$.
Let $E_\rho = C^\omega(A_\rho)\cap C^0(\overline{A_\rho})$ be the
Banach space of analytic functions on the open annulus extending
continuously to its closure.

Let $j: E_\rho \rr \BV$ denote the natural injection. It is continuous
because by a Cauchy estimate, $|f|_\infty + \int_{S^1} |f'| \, dx \leq
(1 + \frac1{\rho}) \|f\|_{E_\rho}$. Write $P^E$ for the restriction of
$P_\delta$ to $E_\rho$.  Then
  \begin{equation}
        j \circ P^E = P_\delta \circ j.
  \end{equation}
The contraction property of $\tau$ shows that
$P^E : E_\rho \rr E_{\th \rho}$ is 
norm-bounded (by $2\th$), and since
the natural injection $E_{\th \rho} \rr E_\rho$ is compact
(in fact nuclear) the operator 
$P^E$ is compact when acting upon $E_\rho$.

Returning now to the above left eigenfunction, $\ell_\delta$ we see that
 \begin{equation}
     0 = \ell_\delta (\lambda_\delta - P_\delta) \circ j
       = \ell_\delta \circ j \circ (\lambda_\delta - P^E).
  \end{equation}
The functional $\ell_\delta \circ j$ is therefore in
the kernel of $\lambda_\delta - (P^E)^*$. On the other hand,
it does not vanish when acting upon the test functions $\SS$
(because $\ell_0$ did not and $\delta$ is small enough)  so that
$\lambda_\delta$ must be in the spectrum of $P^E$ as well.
But $P^E$ is compact so $\lambda_\delta\neq0$ is
necessarily an isolated eigenvalue
of finite multiplicity. Hence, there must be a corresponding
right eigenvector, $\phi^E\in E_\rho$. But then
  \begin{equation}
       (\lambda_\delta-P_\delta) \circ j  (\phi^E)=
        j \circ ((\lambda_\delta-P^E) \phi^E) = 0
  \end{equation}
  shows that $j(\phi^E)$ is indeed the corresponding eigenfunction in
  $\BV$ for the eigenvalue $\lambda_\delta$. And $j(\phi^E)$ is
  manifestly real-analytic on the circle. This finishes the proof of
  Theorem~\ref{theo:main}.

\end{document}